\documentclass{amsart}
\usepackage[utf8]{inputenc}
\usepackage{mathtools,amssymb}
\usepackage{xcolor}
\usepackage{bm}
\newcommand{\be}{\begin{equation}}
\newcommand{\ee}{\end{equation}}
\newcommand{\ba}{\begin{array}}
\newcommand{\ea}{\end{array}}

\newcommand{\bea}{\begin{eqnarray*}}
\newcommand{\eea}{\end{eqnarray*}}
\newcommand{\bean}{\begin{eqnarray}}
\newcommand{\eean}{\end{eqnarray}}

\newcommand{\R}{\bm{R}}

\newcommand{\N}{\bm{N}}

\newcommand{\supp}[1]{{\rm supp}\,#1}

\newcommand{\h}{\mathfrak{H}}
\newcommand{\M}{\mathfrak{M}}
\newcommand{\spr}{\textrm{Spr}}
\newtheorem{theorem}{Theorem}

\newtheorem{proposition}{Proposition}
\newtheorem{lemma}{Lemma}
\newtheorem{remark}{Remark}
\title[Series solutions of integral equations in scattering]{On the series solutions of\\integral equations in scattering}
\author{Mirza Karamehmedovi\'c and Faouzi Triki}

\begin{document}

\begin{abstract} We study the validity of the Neumann or Born series approach in solving the Helmholtz equation 
and coefficient identification  in  related inverse scattering problems.  Precisely, we derive a sufficient and necessary condition under which the series is strongly convergent. We also investigate the rate of convergence of the series. The 
obtained condition is optimal and it can be much weaker than the traditional requirement for the convergence of the series. 
Our approach makes use of reduction space techniques proposed by Suzuki \cite{Suzuki-1976}. Furthermore 
we propose an interpolation method that allows the use of the Neumann series in all cases. Finally, we provide several numerical tests with different medium functions and frequency values  to validate our theoretical results. 

\end{abstract}


\maketitle

\section{Introduction and main results}\label{sec:intro}
Let $d=1,2,\dots$, fix a positive $k_0$ and $\widehat{k}\in S^{d-1}$, let $q\in L^{\infty}(\R^d)$ be compactly supported with $q(x)>-1$, and consider the following Helmholtz equation in $\R^d$ with the Sommerfeld radiation condition:
\begin{equation}\label{eqn:HS}
\left\{\begin{array}{rcl}
(\Delta+k_0^2(1+q(x)))u&=&-k_0^2q(x)e^{ik_0\widehat{k}.x}\quad\text{in}\,\,\,\R^d,\\
\lim_{|x|\rightarrow\infty}|x|^{(d-1)/2}\left(\partial_{|x|}-ik\right)u&=&0\quad\text{uniformly in}\,\,\,x/|x|\in S^{d-1}.
\end{array}\right.
\end{equation}
Convolving the PDE in~\eqref{eqn:HS} with the outgoing fundamental solution\footnote{Here $H_0^{(1)}$ the Hankel function of the first kind and order zero.} of the Helmholtz operator $\Delta+k_0^2$ in $\R^d$,
\[
\Phi_d(x)=\begin{cases}(-2\pi|x|)^{(1-d)/2}(2ik_0)^{-1}\partial_{|x|}^{(d-1)/2}e^{ik_0|x|},\quad&x\in\bm{R}^d\setminus\{0\},\,\,\,d\,\,\,\text{odd},\\(-2\pi|x|)^{(2-d)/2}(4i)^{-1}\partial_{|x|}^{(d-2)/2}H_0^{(1)}(k_0|x|),\quad&x\in\bm{R}^d\setminus\{0\},\,\,\,d\,\,\,\,\text{even},\end{cases}
\]
and integrating by parts, we get the Lippmann-Schwinger equation
\begin{equation}\label{eqn:LSE}
(I-V_q(k_0))u=V_{q}(k_0)e^{ik_0\widehat{k}(\cdot)}\quad\text{in}\,\,\,\R^d,
\end{equation}
where 
\[
V_q(k_0)u(x)=k_0^2\int_{y\in\supp{q}}\Phi_d(x-y)q(y)u(y)dy
\]
exists as an improper integral for each $x\in\R^d$. It is well-known that the integral equation~\eqref{eqn:LSE} is equivalent with~\eqref{eqn:HS}, and that it suffices to solve~\eqref{eqn:LSE} in, say, a 
bounded open ball $B\subset\R^d$ that includes $\supp{q}$, followed by the continuous extension $u(x)=V_q(k_0)[u(\cdot)+\exp(ik_0\widehat{k}(\cdot))](x)$ for $x\in\R^d\setminus B$. The mapping $V_q(k_0):L^2(B)\rightarrow L^2(B)$ is compact, and we shall in the following consider only the restriction of the Lippmann-Schwinger equation in~\eqref{eqn:LSE} to $B$.  
The objective of the paper is to study the successive approximations for  solving the  integral equation \eqref{eqn:LSE}:
\begin{equation}\label{iterative scheme}
u_0 = V_q(k_0)e^{ik_0\widehat{k}}; \;   u_{n+1} = u_0+ V_q(k_0)u_n, \;\; n\in \N. 
\end{equation}
 The computational advantage of  this iterative method  is that it does need  to solve the partial differential equation \eqref{eqn:HS} in the whole space and deal with the radiation conditions. Instead, one can 
 obtain a good approximation $u_n$ of the solution $u$ by applying  successively the integral operator $V_q(k_0)$ if the
 sequence converges.\\

 On the other hand the strong convergence of the sequence $(u_n)_{n\in \N}$ to the solution $u$ of 
the integral equation \eqref{eqn:LSE} is equivalent to the convergence 
of the Neumann series: 
\begin{equation} \label{Neumann series}  \lim_{n\to \infty} u_n =\sum_{j=0}^\infty 
V_q^{j+1}(k_0)e^{ik_0\widehat{k}\cdot} =(I-V_q(k_0))^{-1}V_q(k_0)e^{ik_0\widehat{k}\cdot}.
\end{equation}

In inverse scattering problems the Neumann series approach  known more under the name of {\it Born approximation } was initially employed  to study the quantum mechanical inverse backscattering problem in one dimension (see for instance \cite{MS3} and references therein). The principal advantage 
of using this technique in inverse medium problem is that it
requires solving a linear equation instead
of an oscillatory nonlinear one \cite{BT1, BT2}. It  has  also been applied to various other inverse  problems, including optical  and electrical impedance tomographies, acoustic and electromagnetic parameters identification \cite{BCM,  Kleinman1990, AGKLS, MS, MS2, PMCS, MaS, KMS, KMS2, BG, AMS}.  However, it is important to note, that the strategies considered in these works are based on purely formal analysis or assume strong conditions on the targeted  physical parameters. \\

It is well known that a sufficient condition for the convergence of the Neumann series \eqref{Neumann series} is that the spectral radius $\spr(V_q(k_0))$ of the compact operator $V_q(k_0)$ is strictly less than one, that is $\spr(V_q(k_0))<1$. But this latter condition while it is optimal for the expansion of  the operator $(I-V_q(k_0))^{-1}$ in $L^2(B)$,  it is obviously too restrictive for the convergence of \eqref{Neumann series}.  Then is it possible to derive a necessary and sufficient condition for the convergence of only \eqref{Neumann series}? On the other hand the strong convergence
\begin{equation} \label{optimal condition}
V_q(k_0)^{j}e^{ik_0\widehat{k}\cdot} \rightarrow 0, \; \; j \to \infty,
\end{equation}
in $L^2(B)$, is evidently a necessary condition for the convergence of the series 
\eqref{Neumann series}.  Suzuki in his seminal work  \cite{Suzuki-1976} wondered if this condition is also sufficient. Surprisingly, it turns out that this condition 
also guarantees the convergence of the series. The main idea of the proof is to derive a minimal invariant space
where the expansion of the restriction of $(I-V_q(k_0))^{-1}$ to that space is equivalent to the convergence of the series \eqref{Neumann series}. \\
Let
\begin{equation} \label{ reduced space}
L^2_{k_0,\hat k}(B) = \textrm{Span}\left(V_q(k_0)e^{ik_0\widehat{k}\cdot}, 
V_q(k_0)^{2}e^{ik_0\widehat{k}\cdot} , \cdots, V_q(k_0)^{j+1}e^{ik_0\widehat{k}\cdot} , \cdots \right).
\end{equation}
By construction  $L^2_{k_0,\hat k}(B)$ is invariant by $V_q(k_0)$. Denote $\widetilde V_q(k_0)$ the restriction of 
$V_q(k_0)$ to $L^2_{k_0,\hat k}(B)$. 
Suzuki showed that condition \eqref{optimal condition}  indeed implies $\spr(\widetilde V_q(k_0))<1$, and hence ensures the convergence of the
Neumann series to 
the unique solution.

\begin{proposition}
\label{main result}
The convergence of the Neumann series \eqref{Neumann series} is equivalent to
the condition \eqref{optimal condition}. 
\end{proposition}
 
\begin{remark}\label{rem:1}
Since  $V_q(k_0)$ is a compact operator  the strong convergence  \eqref{optimal condition} can be replaced  by a weak convergence of a subsequence.
Notice that $L^2_{k_0,\hat k}(B)$ can also be generated by finite sums of the sequence  

\[
L^2_{k_0,\hat k}(B) = \textrm{Span}\left(  \sum_{j=0}^J V_q(k_0)^{j+1}e^{ik_0\widehat{k}\cdot};\; \; J\in \N \right).
\]
Recall that the traditional  condition  to ensure the convergence of the Neumann series is \cite{BT2}
\begin{equation} \label{norm condition}
 \|V_q(k_0)\|\leq  C_{k_0,q} = \left(\int_{B}\int_{B}|k_0^2\Phi_d(x-y)q(x)|^2dxdy\right)^{1/2}<1.
\end{equation}
This condition occurs in the  situation for weak scattering, and is not
valid for high wave number $k_0$, or large magnitude of the medium function
$q$. But since  $e^{ik_0\widehat{k}\cdot}$ is sparse 
we expect that $L^2_{k_0,\hat k}(B)$ has a lower dimensionality than the whole space 
$L^2(B)$, and consequently the convergence of the Neumann series \eqref{Neumann series} may happen beyond the conventional 
limitation \eqref{norm condition}. In other words $\spr(\widetilde V_q(k_0))<1$ can be satisfied by a larger  class of wave numbers and 
medium functions not necessary within the weak scattering regime. This was also  observed in many numerical experiments in the past, has
fueled many discussions and  was the origin of several investigations \cite{MS, MS2, PMCS, MaS, KMS, KMS2, AMS}. This pattern
is clearly  confirmed by many numerical examples in section 4.

\end{remark}
\begin{theorem} \label{main2}
   Assume that  the condition \eqref{optimal condition} is satisfied, that is 
   $$ \lim_{n\to \infty}\|V_q(k_0)^{n}e^{ik_0\widehat{k}\cdot}\|=0.$$ Then there exists
   a constant $C>0$ independent of $n$ such that  the following
   error estimate 
   \begin{equation}
       \|u-u_n\| \leq C \|V_q(k_0)^{n}e^{ik_0\widehat{k}\cdot}\|,
   \end{equation}
   holds for all $n \in \N.$
   
\end{theorem}

The rest of the paper is organized as follows. In section 2,
we provide the proofs for  Proposition \ref{main result}, and Theorem \ref{main2}. Section 3 is devoted to the construction of a preconditioner for the integral equation \eqref{eqn:LSE}. Precisely, we propose an interpolation method that allows the use of the Neumann series independently of the fact that the condition \eqref{optimal condition} is fulfilled or not.  We present then several numerical experiments to show the effectiveness of the  derived theoretical results
in section 4.

\section{Proof of the main results} \label{section 2}
In this section we shall prove the main results of the paper. To ease the notation we 
set 
$$
A = V_q(k_0);\;\; \psi = V_q(k_0)e^{ik_0\widehat{k}\cdot}; \;\; \h =L^2(B); \;\;\h_0=
 L^2_{k_0,\hat k}(B).
 $$

\begin{proof} [Proof of Proposition \ref{main result}] If the series 
\begin{equation} \label{abstract Neumann series}
\sum_{j=0}^\infty A^j\Psi = (I-A)^{-1}\Psi,
\end{equation}
is convergent then obviously 
we will have  $A^{j}\Psi\rightarrow0$ strongly in $\h$.
Now assume that 
 $A^{j}\Psi$ converges strongly to zero in $\h$,
and focus on the 
nontrivial opposite direction.  \\

The main observation
of Suzuki is that the convergence of the series \eqref{abstract Neumann series}
depends more on the specific local behavior of the operator $A$ relative to the 
given vector $\psi$ rather than its global properties on the whole space $\h$
which requires that $\spr(A) <1$. Indeed 
consider the  Hilbert subspace $\h_0 \subset \h$ space generated by the vectors $A^j\Psi, \, j\in \N $, that is 
\begin{equation} 
\h_0 =  \textrm{Span}\left(\Psi, A\Psi, \cdots, A^j\Psi, \cdots \right).
\end{equation}
Clearly $\h_0$ is invariant by $A$, and  since $\Psi$ lies in $\h_0$ to prove
that the series \eqref{abstract Neumann series} strongly   in $\h_0$
it is sufficient to show that $A_0$ the restriction of $A$ to $\h_0$ verifies
$\spr(A_0)<1$. Remark that since $\h_0 \subset \h$ we have $\spr(A_0) \leq \spr(A)$. \\

Let $\M$ be the linear 
manifold formed by the  vectors $v \in \h_0$ satisfying $A^jv$ tends  strongly to zero as $j\to \infty$. 
We first remark  that the fact $j\to A^{j}\Psi$ tends strongly to zero,  $\M$ contains  all the vectors $A^j\Psi, \, j\in \N $, and consequently is dense in $\h_0$. \\

Let $\sigma(A_0)$ denotes the spectrum of $A_0$, and set $\sigma_-(A_0) =
\{\lambda \in \Sigma(A_0); \; |\lambda| <1 \}$ and $\sigma_+(A_0) =
\{\lambda \in \Sigma(A_0); \; |\lambda| \geq 1 \}$. Since $A_0$ is compact 
$\sigma_+(A_0)$ is finite, in addition there exists a  rectifiable Jordan curve $\mathcal C_+$  in the resolvent set surrounding  $\sigma_+(A_0)$
and does not contain other eigenvalues of $\sigma(A_0)$. Similarly  there exists  a  rectifiable Jordan curve $\mathcal C_-$   in the resolvent set surrounding  only $\sigma_-(A_0) $. 
Then following \cite{Kato}, the spectral projections 
\begin{eqnarray}
P_\pm = \frac{1}{2i\pi}\int_{\mathcal C_\pm} (\lambda I- A_0)^{-1} d\lambda,
\end{eqnarray}
verify the following identities
\begin{equation} \label{identities}
  P_- +P_+ = I; \; P_-P_+= P_+P_-= 0; \; P_\pm A = A P_\pm.    
\end{equation}

Recalling that $\spr(A_0) = \sup_{\lambda \in \sigma(A_0)}{|\lambda|}.$  Since $\sigma(A_0)$  is a sequence of complex
values that  may converge to zero,  proving that $\spr(A_0)<1$ is equivalent to show that $\sigma_+(A_0)$ is an empty set. Let now 
$v \in P_+ \h_0$. By the density of the set $\M$ in $\h_0$, there 
exists a sequence $(v_n)_{n\in \N_0} \in \M$ that converges 
strongly to $v$. $\N_0 $ here is the set $\N\setminus\{0\}$. Denote $v_{n, \pm} = P_\pm v_n$. Therefore 
$v_n = v_{n, +}+v_{n, -}$. Remarking that $ A v_{n, -}=   A P_- v_n=
P_-Av_n$ converges strongly to  $P_- v =0$ leads to $v_{n, -} \in \M$. 
Hence $ v_{n, +} = v_n - v_{n, -}$  lies in fact in $\M \cap P_+ \h_0$. This shows that  $\M \cap P_+ \h_0$ is indeed dense in $ P_+ \h_0$. But since $\sigma_+(A_0)$ is finite $ P_+ \h_0$ is finite 
dimensional space and consequently $\M \cap P_+ \h_0= P_+ \h_0$.  
This is clear not correct if $P_+ \h_0$ is not trivial (take any
eigenvector of $A_0$ associated to $\lambda \in \sigma_+(A_0)$, 
it obviously does not belong to $\M$). Then $\sigma_+(A_0)$ is an empty set, and finally $\spr(A_0)<1$, which achieves the proof.

\end{proof}
\begin{proof} [Proof of Theorem \ref{main2}] 
Since $\psi \in \M$ we deduce from Proposition \ref{main result} that the Neumann series \eqref{abstract Neumann series} is convergent.
On the other hand we deduce from \eqref{identities} $A^j \psi =  A^j P_-\psi = A_0^j \psi$. Therefore 
\begin{equation} \label{bound}
\|\sum_{j=0}^\infty A^j\Psi\| =  \|\sum_{j=0}^\infty A_0^j\Psi\| = \|(I-A_0)^{-1}\Psi\| \leq  \|(I-A_0)^{-1}\|\|\Psi\|.
\end{equation}
Let $n\in \N$ be fixed. The fact that 
$\psi \in \M$ implies that $A^{n+1}\psi \in \M$. Applying  inequality \eqref{bound} to the 
vector $A^{n+1}\psi$ leads to 
\begin{equation}
\|\sum_{j=n+1}^\infty  A^j\Psi\|  = \|\sum_{j=0}^\infty  A^j A^{n+1}\Psi\| \leq  \|(I-A_0)^{-1}\|\|A^{n+1}\Psi\| \leq
C \|A^{n}\Psi\|,
\end{equation}
with $C = \|(I-A_0)^{-1}\| \|A_0\|$, which finishes the proof of the Theorem.

\end{proof}

\begin{remark}  The proofs stay valid for any general  compact operator $A$ and even  if $\h$ is a  Banach space. 
In the particular case where $A$ is in addition normal, that is $AA^* = A^*A$, the obtained results are straightforward. Indeed
if $\sigma(A)=\{\lambda_k; k\in \N_0\}$ the eigenvalues of $A$, and $P_k$ is the orthogonal spectral  projection 
associated to $\lambda_k$, we have 
\[A = \sum_{k=1}^\infty \lambda_k P_k,\]
and it is clear that the condition $\psi \in \M$  is equivalent to
\[
\psi = \sum_{|\lambda_k|<1}  P_k \psi. 
\]
Therefore
\begin{eqnarray*}
\left\|\sum_{j=n+1}^\infty  A^j\Psi \right\|^2 = \left \| \sum_{|\lambda_k|<1} \frac{\lambda_k^{n+1}}{1-\lambda_k} P_k \Psi \right \|^2 
= \sum_{|\lambda_k|<1} \frac{\lambda_k^{2(n+1)}}{(1-\lambda_k)^2} \|P_k\Psi\|^2\\
\leq \frac{r_0^{2}}{(1-r_0)^2}\|A^{n}\Psi\|^2,
\end{eqnarray*}
where $r_0 = \max_{|\lambda_k|<1}|\lambda_k| = \|A_0\| =\spr(A_0).$ One can verify that $
C = \|(I-A_0)^{-1}\|\|A_0\| = \frac{r_0}{1-r_0}.$
Finally it is easy  to find examples of $A$ such that 
the inequalities
\[
\spr(A_0) \ll 1 \ll \spr(A) =\|A\|,
\]
hold, and where the benefit of considering the reduced space $\h_0$ is indeed remarkable.
\end{remark}

\section{Preconditioning}
By 'preconditioning' we here mean the transformation of the original Lippmann-Schwinger equation $(I-V_q(k_0))u=\psi$ to an integral equation solvable by a convergent Neumann series regardless of the value of $\|V_q(k_0)\|_{L^2(B)\rightarrow L^2(B)}$ and of whether or not the sequence $(\|V_q(k_0)^j\psi\|_{L^2})_{j\in\N}$ converges to zero. See~\cite{Engl-Nashed-1981,Kleinman1988,Kleinman1988-2, Kleinman1990,Osnabrugge2016} for related approaches. Throughout this section we assume the problem dimension $d\in\{1,2,3\}$.
\begin{lemma}\label{lemma:precon}
    If $q(x)\ge0$, $q\not\equiv0$, then there is a complex constant $\gamma$ such that the solution of the equation $(I-V_q(k_0))u=\psi$ in $L^2(B)$ is expressible in terms of the convergent Neumann series
    \[
    u=\sum_{j=0}^{\infty}M^j\gamma\psi,
    \]
    where $M=(1-\gamma)I+\gamma V_q(k_0)$.
\end{lemma}
\begin{proof}
    Let $V_q(k_0)\varphi=\lambda\varphi$ in $\bm{R}^d$ with nonzero $\lambda$ and with $\varphi$ not identically zero outside any bounded ball in $\R^d$. Then
\begin{equation}
\left\{\begin{array}{rcl}
\varphi''+k_0^2(1+q(x)/\lambda)\varphi&=&0,\qquad\quad x\in]-R,R[,\\
-\varphi'(-R)&=&ik_0\varphi(-R),\\
\varphi'(R)&=&ik_0\varphi(R)
\end{array}\right.
\end{equation}
for $d=1$, and
\begin{equation*}
\left\{\begin{array}{rcl}
(\Delta+k_0^2(1+q(x)/\lambda))\varphi&=&0\quad\text{in}\,\,\,\R^d,\\
\lim_{|x|\rightarrow\infty}|x|^{(d-1)/2}\left(\partial_{|x|}-ik\right)\varphi&=&0\quad\text{uniformly in}\,\,\,x/|x|\in S^{d-1},
\end{array}\right.
\end{equation*}
for $d\in\{2,3\}$. Thus, for sufficiently large $R>0$ we have
\begin{align}\label{eqn:former}
0&=\int_{|x|<R}(\overline \varphi\Delta \varphi+k_0^2|\varphi|^2+k_0^2\lambda^{-1}q(x)|\varphi|^2)\nonumber\\&=\int_{|x|<R}(k_0^2|\varphi|^2-|\nabla \varphi|^2)+\int_{|x|=R}\overline{\varphi}\partial_r\varphi+k_0^2\lambda^{-1}\int_{x\in\supp q}q(x)|\varphi|^2,
\end{align}
as well as $(\Delta+k_0^2)\varphi=0$ in $\{|x|>R\}$. In the case $d=1$ we readily find that
\[
\Im\int_{|x|=R}\overline{\varphi}\partial_r\varphi=k_0(|\varphi(-R)|^2+|\varphi(R)|^2)>0,
\]
while in the case $d\in\{2,3\}$ we can follow the argument in, e.g.,~\cite[Theorem 2.13, p. 38]{CK4} to find
\[
\Im\int_{|x|=R}\overline{\varphi}\partial_r\varphi>0.
\]
Hence 
\[
\int_{x\in\supp q}q(x)|\varphi|^2dx >0, 
\]
and this in conjunction with~\eqref{eqn:former} gives
\[
\Im(\lambda^{-1})=-\frac{\Im\int_{|x|=R}\overline{\varphi}\partial_r\varphi}{k_0^2\int_{x\in\supp q}q(x)|\varphi|^2}<0,
\]
so $\Im\lambda>0$ and finally $\Re(e^{i\alpha}(1-\lambda))>0$ if $\alpha\in]-\pi/2,\pi/2[$ satisfies
\begin{equation}\label{eqn:condition}
\tan\alpha>\max_{\lambda'\in\sigma(V_q(k_0))}\frac{\Re\lambda'-1}{\Im\lambda'},
\end{equation}
where $\sigma(V_q(k_0))$ is the spectrum of $V_q(k_0)$. The existence of the maximum in~\eqref{eqn:condition} follows from the fact that the eigenvalues of the compact operator $V_q(k_0):L^2(B)\rightarrow L^2(B)$ reside in the closed ball $\{\lambda'\in\bm{C},\,\,\,|\lambda'|\le\|V_q(k_0)\|_{L^2(B)\rightarrow L^2(B)}\}$ and can accumulate only at zero. These facts also imply that
there exists $\varepsilon>0$ such that $|\gamma(1-\lambda')-1|<1$ for all $\lambda'\in\sigma(V_q(k_0))$, where $\gamma=\varepsilon e^{i\alpha}$. It remains to notice that 
each eigenvalue $\mu$ of $M$ is of the form $\mu=1-\gamma+\gamma\lambda'$, where $\lambda'$ is some eigenvalue of $V_q(k_0)$, and finally that the equation $(I-M)u=\gamma\psi$ is equivalent with the equation $(I-V_q(k_0))u=\psi$.
\end{proof}

We can in fact be more specific in a special case in dimension one. Let $L$ be a positive constant, set $B=]0,L[$, and let $q(x)\equiv q_0={\rm const.}>0$ for $x\in\overline{B}$.
\begin{lemma}\label{lemma:precon2}
If $\varepsilon'>k_0Lq_0/2$,
\[
\alpha=\arctan\frac{1+\frac{k_0Lq_0}{2}\varepsilon'}{\varepsilon'-\frac{k_0Lq_0}{2}},
\]
\[
0<\varepsilon<\frac{1}{{1+k_0 Lq_0/2}}\frac{|\tan(2\max\{\alpha,\arctan\varepsilon'\})|}{\tan(\max\{\alpha,\arctan\varepsilon'\})}
,
\]
and $\gamma=\varepsilon e^{i\alpha}$,
then each eigenvalue $\mu$ of $M=(1-\gamma)I+\gamma V_q(k_0)$ satisfies $|\mu|<1$.
\end{lemma}
\begin{proof}





Assume $\lambda\in\bm{C}\setminus\{0\}$ and $\varphi\in L^2(]0,L[)$, $\varphi\not\equiv0$, satisfy 
\begin{equation}\label{eqn:sys}
V_q(k_0)\varphi(x)=\lambda\varphi(x),\quad x\in(0,L).
\end{equation}
Integration by parts readily shows the equivalence of the Lippmann-Schwinger equation~\eqref{eqn:sys} with the Helmholtz system
\begin{equation}\label{eqn:hs}
\left\{\begin{array}{rcl}
\varphi''(x)+k_0^2s^2\varphi(x)&=&0,\qquad\quad x\in]0,L[,\\
-\varphi'(0)&=&ik_0\varphi(0),\\
\varphi'(L)&=&ik_0\varphi(L),
\end{array}\right.
\end{equation}
where $s=\left(1+q_0/\lambda\right)^{1/2}$.
The eigenvectors of the Laplacian on $]0, L[$ are generally of the form
\begin{equation}\label{eqn:form}
\varphi(x)=A\exp(ik_0sx)+B\exp(-ik_0sx),
\end{equation}
and we then readily find that~\eqref{eqn:hs} is equivalent with~\eqref{eqn:form} together with \begin{equation}\label{eqn:equiv}
s\neq1,\quad B=A(s+1)/(s-1),\quad e^{2ik_0 Ls}=(s+1)^2/(s-1)^2.
\end{equation}
Next define the constant $T(k_0 L)>0$ by $T(k_0 L)\sinh(k_0 LT(k_0 L))=1$. Using the last condition in~\eqref{eqn:equiv}, 
we find that, necessarily,
\[
s\in S_{\pm}(k_0 L)=\left\{\frac{-\cosh k_0 Lt\pm\sqrt{1-t^2\sinh^2k_0 Lt}}{\sinh k_0 Lt}+it,\,\,t\in(0,T(k_0 L)]\right\},
\]
which in turn implies
\[
\lambda\in\Lambda(q_0,k_0 L)=\left\{\frac{q_0}{s^2-1},\,\,s\in S_-(k_0 L)\cup S_+(k_0 L)\right\}.
\]
As an example, Figure~\ref{fig:example_spectrum} shows the  set $\Lambda(q_0,k_0L)$, as well as the numerically computed spectrum, for two sets of parameter values for $k_0$, $L$, and $q_0$.
\begin{figure}
    \centering
    \includegraphics[scale=0.3]{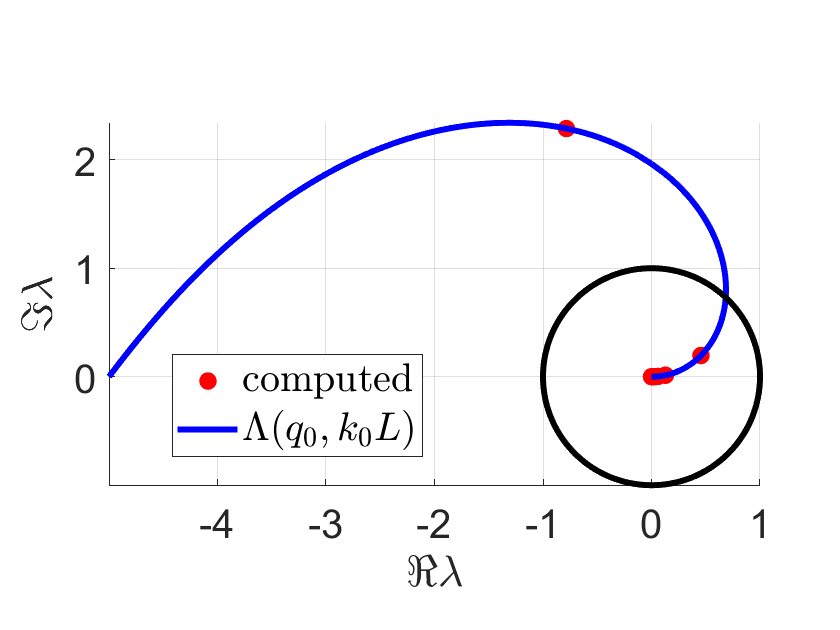}
    \includegraphics[scale=0.3]{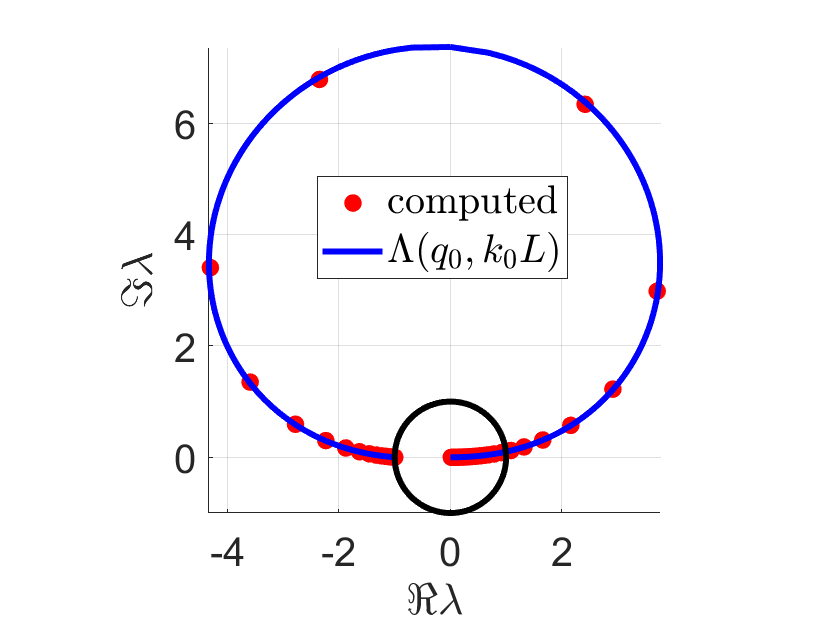}
    \caption{The theoretically predicted set $\Lambda(q_0,k_0L)$ that includes the eigenvalues of $V_q(k_0)$, plotted against the numerically computed eigenvalues. The parameter values are: top, $k_0=1$, $L=1$, and $q_0=5$; bottom, $k_0=50$, $L=1$, and $q_0=1$. The reference circle is centered at the origin and has radius one.}
    \label{fig:example_spectrum}
\end{figure}

Since $\Im\lambda>0$ and since the eigenvalues of $V_q(k_0)$ accumulate precisely at zero, we have
\[
\lim_{t\rightarrow0}\frac{\Re\omega(t)-1}{\Im\omega(t)}=-\infty
\]
for
\[
\omega(t)=\frac{q_0}{s(t)^2-1},\quad t\in(0,T(k_0 L)];
\]
here $s(t)$ is given by the above definition of $S_{\pm}(k_0L)$. Furthermore, if $\Re\omega(t)<1$ then $(\Re\omega(t)-1)/\Im\omega(t)<0$, while $\Re\omega(t)\ge1$ implies
\[
\frac{1-t^2\sinh^2k_0 Lt+\cosh(k_0 Lt)\sqrt{1-t^2\sinh^2k_0 Lt}}{\sinh^2k_0 Lt}\le\frac{q_0}{2};
\]
the latter can be seen by rewriting $\Re\omega-1\ge0$ as
\[
((\Re s)^2-(\Im s)^2-1)((\Re s)^2-(\Im s)^2-1-q_0)+4(\Re s)^2(\Im s)^2\le0,
\]
and noting that $(\Re s)^2-(\Im s)^2-1>0$ and $(\Re s)^2(\Im s)^2>0$. Now since also
\begin{gather*}
\frac{1-t^2\sinh(k_0Lt)^2-\cosh(k_0Lt)\sqrt{1-t^2\sinh(k_0Lt)^2}}{t\sinh(k_0Lt)(\cosh(k_0Lt)-\sqrt{1-t^2\sinh(k_0Lt)^2})}\\\le\frac{1-t^2\sinh(k_0Lt)^2+\cosh(k_0Lt)\sqrt{1-t^2\sinh(k_0Lt)^2}}{t\sinh(k_0Lt)(\cosh(k_0Lt)+\sqrt{1-t^2\sinh(k_0Lt)^2})},
\end{gather*}
we have
\[
\frac{\Re\omega(t)-1}{\Im\omega(t)}<\frac{\Re\omega(t)}{\Im\omega(t)}<\frac{1-t^2\sinh^2k_0 Lt+\cosh(k_0 Lt)\sqrt{1-t^2\sinh^2k_0 Lt}}{t\sinh(k_0 Lt)(\cosh k_0Lt+\sqrt{1-t^2\sinh^2k_0 Lt})},
\]
so $\Re\omega(t)\ge1$ implies
\begin{align*}
\frac{\Re\omega(t)-1}{\Im\omega(t)}&<\frac{q_0}{2}\frac{\sinh k_0 Lt}{t(\cosh k_0 Lt+\sqrt{1-t^2\sinh^2k_0 Lt})}\le\frac{q_0}{2}\frac{\tanh k_0 Lt}{t}\\&\le\frac{q_0}{2}\lim_{\tau\searrow0}\frac{\tanh k_0 L\tau}{\tau}=\frac{q_0k_0 L}{2},
\end{align*}
that is, we get a similar estimate on $(\Re\omega(t)-1)/\Im\omega(t)$ as we do on $\|V_q(k_0)\|$. Next, define
\begin{equation}\label{eqn:ang1}
\xi_+=\arctan\sup_{\omega\in\Lambda(q_0,k_0 L)}\frac{\Im\left(e^{i\alpha}(1-\omega)\right)}{\Re\left(e^{i\alpha}(1-\omega)\right)}
\end{equation}
and
\begin{equation}\label{eqn:ang2}
\xi_-=\arctan\inf_{\omega\in\Lambda(q_0,k_0 L)}\frac{\Im\left(e^{i\alpha}(1-\omega)\right)}{\Re\left(e^{i\alpha}(1-\omega)\right)}.
\end{equation}
With $\alpha\in(0,\pi/2)$, $\varepsilon>0$, and $\gamma=\varepsilon e^{i{\alpha}}$, we have $\arg(\gamma(1-\omega))\in[\xi_-,\xi_+]$ for all $\omega\in\Lambda(q_0,k_0 L)$, and if
\begin{equation}\label{eqn:epsilon}
\varepsilon<\frac{1}{{1+k_0 Lq_0/2}}\frac{|\tan(2\max\{|\xi_+|,|\xi_-|\})|}{\tan(\max\{|\xi_+|,|\xi_-|\})}
\end{equation}
then $|\gamma(1-\omega)-1|<1$ for all $\omega\in\Lambda(q_0,k_0 L)$, and specifically $|\gamma(1-\lambda')-1|<1$ for all eigenvalues $\lambda'$ of $V_q(k_0)$; the condition~\eqref{eqn:epsilon} can be deduced by requiring $\Im(\gamma(1-\omega))<|\tan(\pi-2\arg(\gamma(1-\omega)))|$ and using $|1-\omega|\le1+\|V_q(k_0)\|$. Note that we must choose $\alpha<\pi/2$ rather than $\alpha=\pi/2$ since
\[
\frac{\Im(e^{i\alpha}(1-\omega(t)))}{\Re(e^{i\alpha}(1-\omega(t)))}\underset{t\rightarrow0}{\longrightarrow}\tan\alpha,
\]
which for $\alpha=\pi/2$ forces $\xi_+=\pi/2$ and thus $\varepsilon<0$. Note furthermore that we can bound~\eqref{eqn:ang1}--\eqref{eqn:ang2} analytically for
\[
\alpha=\arctan\frac{1+\frac{k_0Lq_0}{2}\varepsilon'}{\varepsilon'-\frac{k_0Lq_0}{2}},
\]
with $\varepsilon'>k_0Lq_0/2$, since
\[
-\cot\alpha\le\frac{\Im\left(e^{i\alpha}(1-\omega)\right)}{\Re\left(e^{i\alpha}(1-\omega)\right)}\le\tan\alpha,\quad\omega\in\Lambda(q_0,k_0 L),\quad\Re\omega\le1,
\]
since furthermore (recall that $(\Re\omega-1)/\Im\omega<k_0Lq_0/2<\tan\alpha$)
\[
-\frac{\frac{k_0Lq_0}{2}\tan\alpha+1}{\tan\alpha-\frac{k_0Lq_0}{2}}\le\frac{\Im\left(e^{i\alpha}(1-\omega)\right)}{\Re\left(e^{i\alpha}(1-\omega)\right)}<0,\quad\omega\in\Lambda(q_0,k_0 L),\quad\Re\omega>1,
\] 
and finally since
\[
-\cot\alpha>-\frac{\frac{k_0Lq_0}{2}\tan\alpha+1}{\tan\alpha-\frac{k_0Lq_0}{2}}=-\varepsilon'.
\]
\end{proof}
\begin{remark}\label{rem:drawback}
A drawback of preconditioning is that, with increasing $k_0Lq_0$, the acceptable values of $\alpha$ and of $\arctan\varepsilon'$ tend to $\pi/2$, and $\varepsilon$ therefore tends to zero. Thus, while the Neumann series remains convergent, the equation
\[
(I-M)u=\gamma\psi
\]
may be said, especially in a numerical context, to lose information about the operator $V_q(k_0)$ and about the original inhomogeneity $\psi$, as both are multiplied with $\gamma=\varepsilon e^{\i\alpha}$ there.
\end{remark}
\begin{remark}\label{rem:better}
Instead of using the bound on $\varepsilon$ stated in Lemma~\ref{lemma:precon2}, we can estimate $\xi_+$ and $\xi_-$ from~\eqref{eqn:ang1}--\eqref{eqn:ang2} numerically and arrive at a larger sufficiently small $\varepsilon$ using~\eqref{eqn:epsilon}.
\end{remark}
We show numerical examples of the use of preconditioning in Section~\ref{sec:preconnum}.

\section{Numerical examples} \label{sec:preconnum}
We here present several numerical examples in dimension one. Fix a positive wavenumber $k_0$, obstacle size $L>0$, medium function $q\in L^2(]0,L[)$, $q(x)>-1$,
and consider the following system for the scattered wave $u(x)$ corresponding to the left excitation $\exp(ik_0x)$ in dimension one:
\begin{equation}\label{eqn:psi1}
\left\{\begin{array}{rcl}
\psi^{\prime \prime}(x)+k_0^2(1+q(x))\psi(x) &=&  -k_0^2q(x)\exp(ik_0x),\quad x\in]0,L[,\\
-\psi^\prime(0)&=&ik_0\psi(0),\\\psi^\prime(L)&=&ik_0\psi(L).
\end{array}\right.
\end{equation}
The function $G(x,y)=(i/2k_0)\exp(ik_0|x-y|)$, $x,y\in[0,L]$, is the free-space Green's function associated with the boundary problem~\eqref{eqn:psi1}, since
\begin{equation*}
\left\{\begin{array}{rcl}
(\partial_y^2+k_0^2)G(x,y)&=&-\delta(x-y),\quad x,y\in]0,L[,\\-\partial_yG(x,0)&=&ik_0G(x,0),\quad x\in]0,L[,\\\partial_yG(x,L)&=&ik_0G(x,L),\quad x\in[0,L].\end{array}\right.
\end{equation*}
Multiplying the differential equation in~\eqref{eqn:psi1} with $G(x,y)$ 
and integrating by parts, we get the Lippmann-Schwinger equation
\begin{equation}
(I-V_q(k_0))\psi(x)=V_q(k_0)\exp(ik_0\cdot)(x),\quad x\in]0,L[,
\end{equation}
where
\[
V_q(k_0)u(x)=\frac{ik_0}{2}\int_{y=0}^Le^{ik_0|x-y|}q(y)u(y)dy,\quad x\in]0,L[.
\]
The operator $V_q(k_0):L^2(]0,L[)\rightarrow L^2(]0,L[)$ is compact, with norm satisfying
\[
\|V_q(k_0)\|^2\le\int_{x=0}^L\int_{y=0}^L\left|\frac{ik_0}{2}e^{ik_0|x-y|}q(y)\right|^2dydx=\frac{k_0^2L\|q\|^2_2}{4}.
\]

\begin{figure}
    \centering
    \includegraphics[scale=0.34]{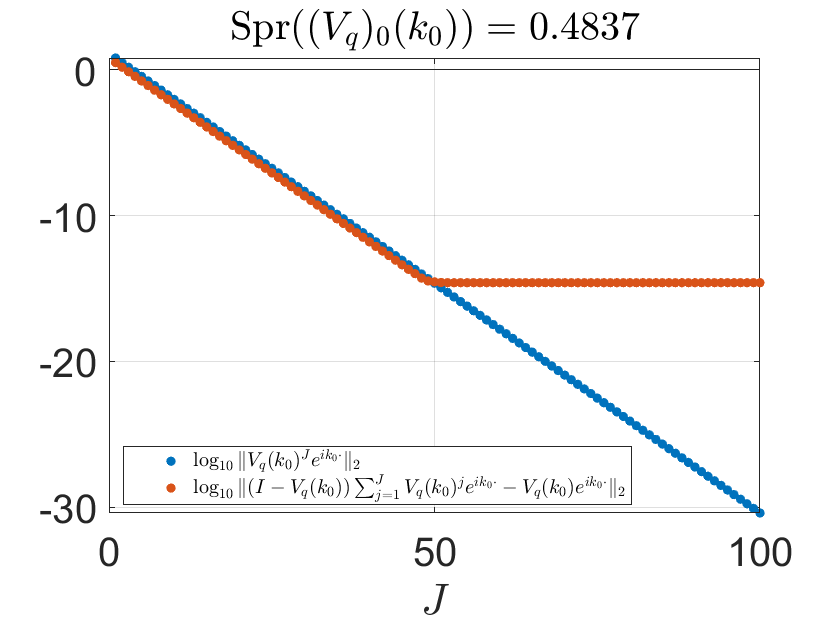}
    \includegraphics[scale=0.34]{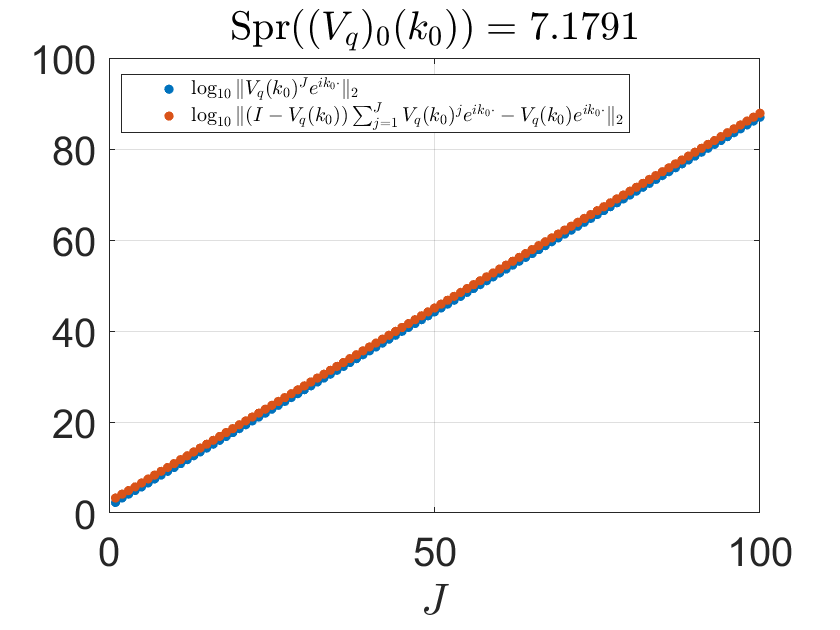}
    \caption{Top: $k_0=1$, $L=1$, $q_0=1$. Bottom: $k_0=50$, $L=1$, $q_0=1$.}
    \label{fig:konverg}
\end{figure}

Now let $q_0 \in ]-1, +\infty[$, and  assume that $q|_{]0, L[}\equiv q_0$. One can easily prove \[\mathfrak H_0= \textrm{Span}\left(V_{q}(k_0)^j e^{ik_0 \cdot}; 
\; j\in \mathbb N\right)  = \textrm{Span} \left(V_{1}(k_0)^j e^{ik_0 \cdot}; \; j\in \mathbb N\right).\]
On the other hand, since $V_q(k_0) = q_0 V_1(k_0)$, we have
\[\textrm{Spr}\left((V_{q})_0(k_0)\right) = 
|q_0|\textrm{Spr}\left((V_1)_0(k_0)\right),\]
where $ (V_{q})_0(k_0)$ and $(V_{1})_0(k_0)$ are the restrictions of $ V_{q}(k_0)$ and $V_{1}(k_0)$, respectively, to the 
space $\mathfrak H_0$. Following the proof of Proposition~\ref{main result},  the Born series 
\[
\sum_{j=1}^\infty \left(V_q(k_0)\right)^j e^{ik_0\cdot} 
\]
is convergent if and only if  $\textrm{Spr}\left((V_{q})_0(k_0)\right)<1$, that is, if and only if
\[
|q_0|<\frac{1}{\textrm{Spr}\left((V_1)_0(k_0)\right)}.
\]
To illustrate this, we show in Figure~\ref{fig:konverg} two cases of repeated application of $V_q(k_0)$ on the original right-hand side $V_q(k_0)e^{ik_0\cdot}$.
In the case where Spr$((V_q)_0(k_0))<1$, the series $(V_q(k_0)^Je^{ik_0\cdot})$ converges strongly to zero, while $(I-V_q(k_0))\sum_{j=1}^JV_q(k_0)^je^{ik_0\cdot}$ converges strongly to $V_q(k_0)e^{ik_0\cdot}$, The second of the two convergence processes plateaus for large values of $J$ due to the effect of numerical errors. In contrast, neither sequence converges in the case where Spr$((V_q)_0(k_0))\ge1$.

\begin{figure}[]
    \centering
    \includegraphics[scale=0.25]{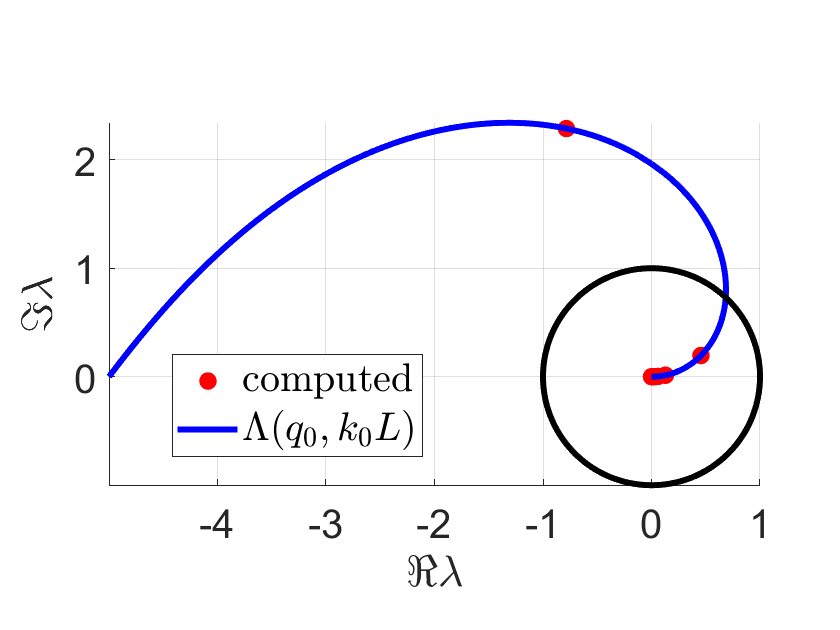}
    \includegraphics[scale=0.25]{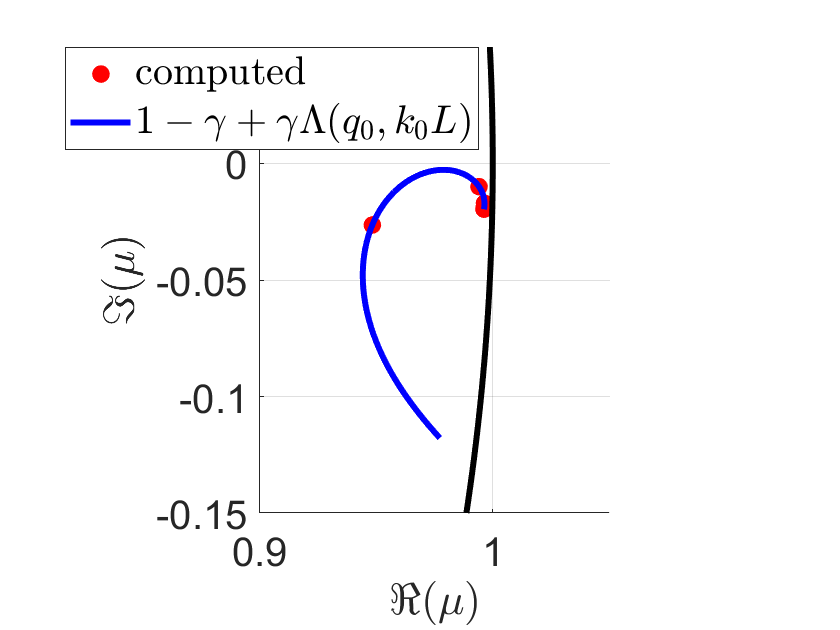}

    \includegraphics[scale=0.25]{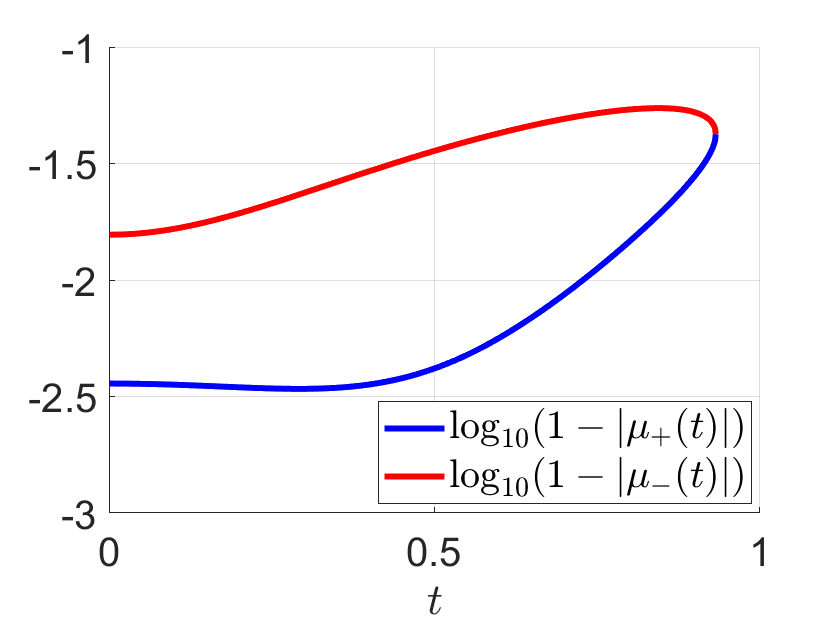}

    a) $k_0$=1, $L=1$, $q_0=5$, $\|V_q(k_0)\|=2.4183$, $\varepsilon'=5.2$, $\alpha=1.3803$, $\varepsilon=0.02$, $T=0.9320$.
    
    \includegraphics[scale=0.25]{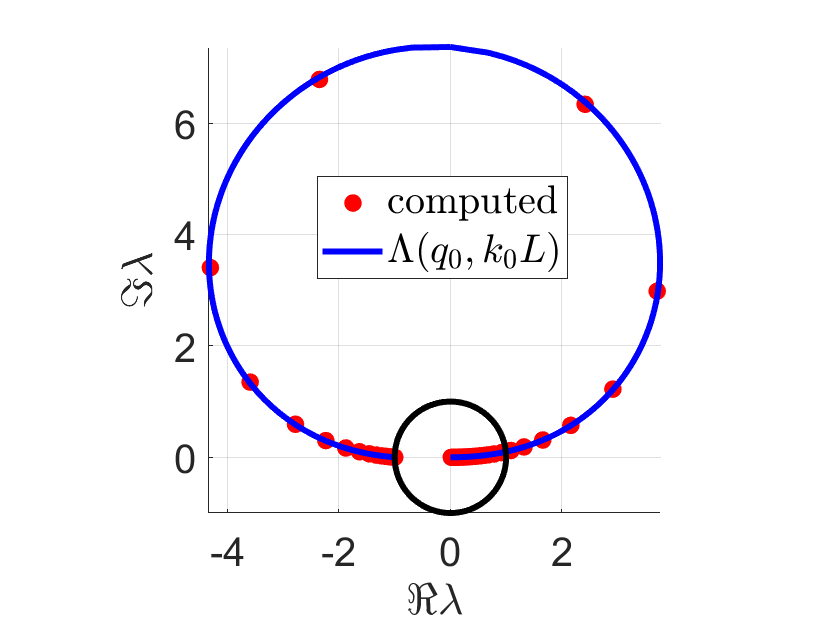}
    \includegraphics[scale=0.25]{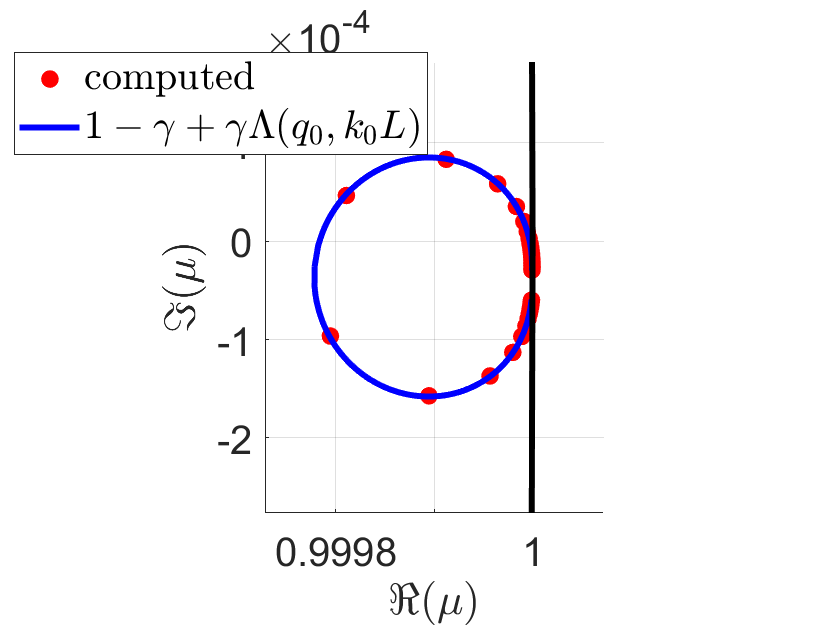}

    \includegraphics[scale=0.25]{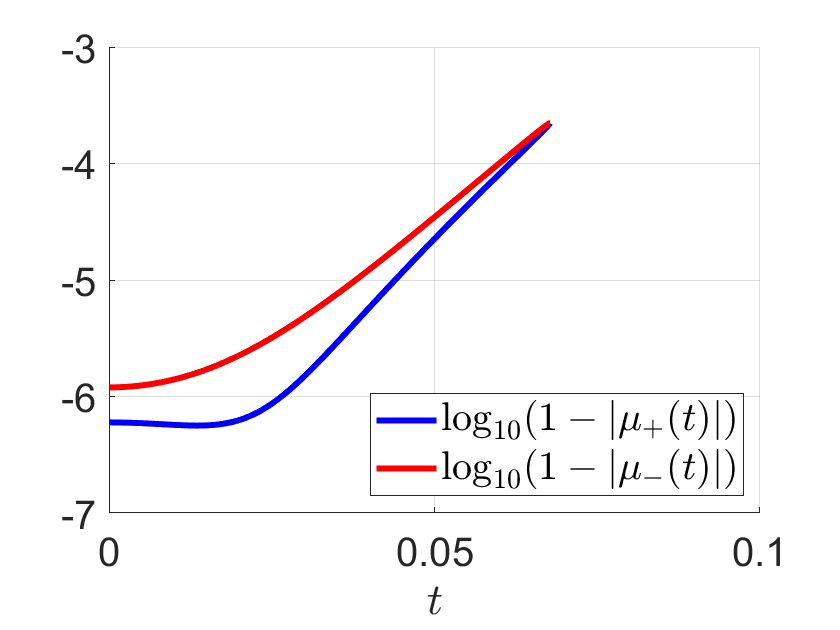}

    b) $k_0=50$, $L=1$, $q_0=1$, $\|V_q(k_0)\|=7.1791$, $\varepsilon'=50$, $\alpha=1.5508$, $\varepsilon=3\cdot10^{-5}$, $T=0.0677$.
    
    \caption{The original theoretically predicted and numerically computed spectra of $V_q(k_0)$; the transformed curve $\{z=1-\gamma(1-\omega),\,\,\,\omega\in\Lambda(q_0,k_0 L)\}\subset\bm{C}$ is included in the unit open disk centered at 1.}
    \label{fig:orig_transf_1}
\end{figure}

\begin{figure}
    \centering
    \includegraphics[scale=0.4]{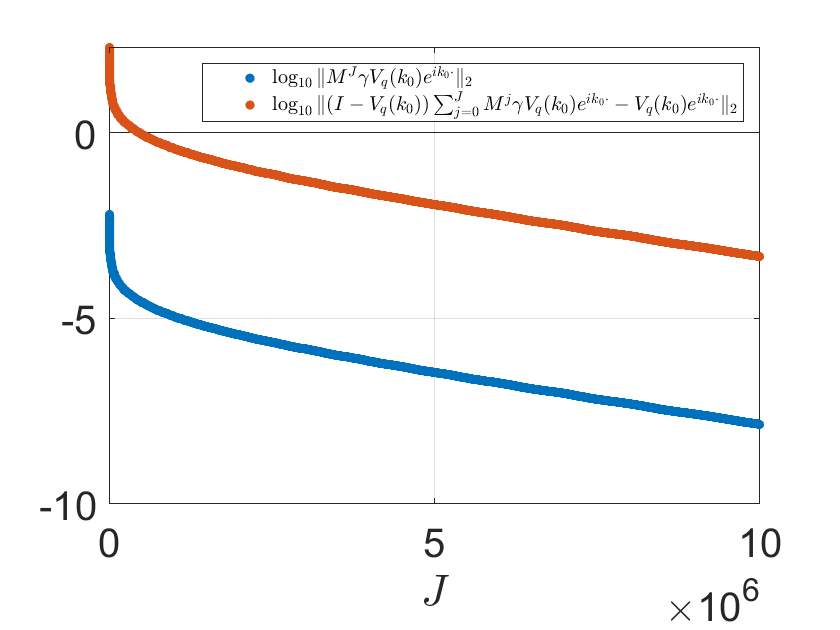}
    \caption{Preconditioning of the equation $(I-V_q(k_0))u=V_q(k_0)e^{ik_0\cdot}$ results in a convergent Neumann series solution. The parameters here are as in Figure~\ref{fig:konverg} (bottom) and Figure~\ref{fig:orig_transf_1} b).}
    \label{fig:precon_works}
\end{figure}

Finally, we illustrate Lemma~\ref{lemma:precon2} numerically in Figures~\ref{fig:orig_transf_1} and~\ref{fig:precon_works}. We indeed get a convergent Neumann series solution for the case $k_0=50$, $L=1$, $q_0=1$, albeit the convergence is rather slow.


\end{document}